\documentclass[12 pt,onesite, openany]{article} 
\usepackage{amsmath,amsxtra,amssymb,latexsym, amscd,amsfonts, indentfirst}
\usepackage[mathscr]{eucal}

\begin{document}
\parskip 3pt
\baselineskip 16pt
\pagenumbering {arabic}
\setcounter{page}{1}
\date{}
\title {\small \bf MIXED MULTIPLICITIES OF MULTI-GRADED ALGEBRAS\\ OVER NOETHERIAN LOCAL RINGS }
\author {\small  Duong Quoc Viet and Truong Thi Hong Thanh\\
\small Department of Mathematics, Hanoi University of Education\\ 
 \small 136 Xuan Thuy Street, Hanoi, Vietnam\\ 
 \small \it E-mail: duongquocviet@fmail.vnn.vn\\
 \small \it E-mail: thanhtth@gmail.com}
\maketitle
\parbox[c]{13 cm}{{\small ABSTRACT: Let $S= \bigoplus_{n_1, ..., n_s \geq 0}S_{(n_1, ..., n_s)}$ be a finitely generated standard multi-graded algebra over a Noetherian local ring $A.$ This paper first expresses mixed multiplicities of $S$ in term of Hilbert-Samuel multiplicity that explained the mixed multiplicities $S$ as the Hilbert-Samuel multiplicities for quotient modules of $S_{(n_1, ..., n_s)}$. As an application, we get formulas for  the mixed multiplicities of ideals that covers the main result of Trung-Verma in [TV].}} 
 \vskip 0.2 cm
\noindent {\bf 1. Introduction}
 
\noindent Throughout this paper, let $(A, \frak{m})$ denote a Noetherian local ring with maximal ideal 
$\mathfrak{m}$, infinite residue $k = A/\mathfrak{m}$; 
$S= \bigoplus_{n_1, ..., n_s \geq 0}S_{(n_1, ..., n_s)}$ 
($s>0$) a finitely generated standard $s$-graded algebra over $A$. 
Let $J$ be an $\mathfrak{m}$-primary ideal of $A$. Set
$$D_J(S) = \bigoplus_{n \geq  0} \frac{J^nS_{(n, ..., n)}}{J^{n + 1}S_{(n, ..., n)}}$$
and $\ell = \dim{D_\mathfrak{m}(S)}$. Then
$$\ell_A(\frac{J^{n_0}S_{(n_1, \ldots, n_s)}}{J^{n_0 + 1}S_{(n_1, \ldots, n_s)}})$$
is a polynomial of total degree $\ell - 1$ in $n_0, n_1, \ldots, n_s$ for all large $n_0, n_1, \ldots, n_s$ (see Section 3). 
If we write the term of total degree $\ell - 1$ in this polynomial in the form
$$\sum_{k_0 + k_1 + \cdots + k_s = \ell - 1} e(J, k_0, k_1, \ldots, k_s, S) 
\frac{n_0^{k_0}n_1^{k_1}\cdots n_s^{k_s}}{n_0!n_1! \cdots n_s!}$$
then $e(J, k_0, k_1, \ldots, k_s, S)$ are non-negative integers not all zero and called the 
{\it mixed multiplicity of $S$ of type $(k_0, k_1, \ldots, k_s)$ with respect to $J$}.
 
In particular, when $S = A[I_1t_1, \ldots, I_st_s]$ is a multi-graded Rees algebra of ideals $I_1, \ldots, I_s$,  $e(J, k_0, k_1, \ldots, k_s, S)$ exactly is the mixed multiplicity of a set of ideals $(J, I_1, \ldots, I_s)$ in local ring $A$(see [Ve2] or [HHRT]).
 
Mixed multiplicities of ideals were first introduced by Teissier and Risler
 in 1973 for two $\mathfrak{m}$-primary ideals and in this case they can be interpreted as the 
multiplicity of  general elements [Te]. Next, Rees in 1984 proved that each mixed multiplicities of a set of $\mathfrak{m}$-primary ideals
 is the multiplicity of a joint reduction of them [Re]. In general, mixed multiplicities have been mentioned in the works of Verma, Katz, Swanson and other authors, see e.g. [Ve1], [Ve2], [Ve3], [Sw], [HHRT], [KV], [Tr1].  By using the concept of (FC)-sequences,
 Viet in 2000 showed that one can transmute mixed multiplicities of a set of arbitrary ideals into
 Hilbert-Samuel multiplicities [Vi]. Trung in 2001 gave the criteria for the positivity of mixed multiplicities of an ideal $I$ [Tr2]. Similar to the methods of Viet [Vi], Trung and Verma in 2007 characterize also mixed multiplicities of a set of ideals, in term of superficial sequences [TV]. Moreover, some another authors have extended mixed multiplicities of a set of ideals to modules, e.g. Kirby and Rees in [KR1, KR2]. Kleiman and Thorup in [KT1, KT2], Manh and Viet in [MV1]. In a recent paper   [VM]  Viet and Manh investigated  the mixed multiplicities of multigraded algebas over  Artinian local rings.
 
   In this paper, we consider  mixed multiplicities of multi-graded algebra $S$ over Noetherian local ring. Our aim is  to characterize mixed  multi-graded 
of $S$ with respect to $J$ in term of Hilbert-Samuel multiplicity (Theorem 3.3, Sect.3). As an application, we get a version of Theorem 3.3 for mixed multiplicities of arbitrary ideals in local rings (Theorem 4.3, Sect.4) that covers the main result in [TV].    
 
The paper is divided in four sections. In Section 2, we investigate (FC)-sequences of multi-graded algebras. Section 3 gives some results on expressing mixed multiplicities of multi-graded algebras in terms of  Hilbert-Samuel multiplicity.  Section 4 devoted to the discussion of  mixed multiplicities of arbitrary ideals in local rings.

  \vskip 0.2cm
\noindent {\bf  2. (FC)-sequences of multi-graded algebras}

The author in [Vi] built (FC)-sequences of ideals in local ring for calculating mixed multiplicities of set of ideals. In order to study mixed multiplicities of multi-graded algebras, this section introduces (FC)-sequences in multi-graded algebras and gives some important properties of these sequences.

Set $\mathfrak{a}: \mathfrak{b}^\infty = \bigcup_{n=0}^\infty (\mathfrak{a} : \mathfrak{b}^n)$, and
\begin{align*}
&(M : N)_A = \{a \in A \mid aN \subset M\};\\
&S_+ = \bigoplus_{n_1 + \cdots  + n_s > 0}S_{(n_1, ..., n_s)} ; \\
&S_i = S_{(0, \ldots, \underset{(i)}{1} , \ldots, 0)};\\
&S_i^+  = S_iS= \bigoplus_{n_i > 0}S_{(n_1, ..., n_s)} (i = 1, 2, ..., s);\\
&S_{++}  = S_1^+ \cap \cdots \cap S_s^+ = \bigoplus_{n_1, \ldots, n_s> 0}S_{(n_1, ..., n_s)} = S_{(1, \ldots, 1)}S.
\end{align*}
\noindent {\bf Definition 2.1.} {\it Let $S= \bigoplus_{n_1, \ldots, n_s \geq 0}S_{(n_1, \ldots, n_s)}$ be a finitely generated standard $s$-graded algebra over a Noetherian local ring $A$ such that $S_{++}$ is non-nilpotent  and let $I$ be an  ideal of $A$. 
A homogeneous element $x \in S$ is called a weak-(FC)-element of $S$ with respect to $I$ 
if there exists $i \in \{1, 2, \ldots, s\}$ such that $x \in S_i$ and
\begin{list}{}{ \setlength{\leftmargin}{1.8cm}\setlength{\labelwidth}{1.3cm} }
\item [\rm (FC$_1$): ] 
$xS_{(n_1, \ldots, n_i - 1, \ldots, n_s)}\cap I^{n_0}S_{(n_1, \ldots, n_s)} 
= xI^{n_0}S_{(n_1, \ldots,n_i - 1, \ldots, n_s)}$
for all large \\$n_0, n_1, \ldots, n_s$. 
\item [\rm (FC$_2$): ] $x$ is a filter-regular element of $S$, i.e.,
$(0 : x)_{(n_1, \ldots, n_s)} = 0$
for all large $n_1, \ldots, n_s$.
\end{list} 
Let $x_1, \ldots, x_t$ be a sequence in $S$. We call that    $x_1, \ldots, x_t$ is  a weak-(FC)-sequence  of $S$ with respect to $I$ if $\bar{x}_{i + 1}$ is a weak-(FC)-element  of $S/(x_1, \ldots, x_{i})S$ with respect to $I$ for all $i = 0, 1, \ldots, t - 1$, where $\bar{x}_{i + 1}$ is the image of $x_{i + 1}$ in $S/(x_1, \ldots, x_{i})S$.
}
 \vskip 0.2 cm
\noindent {\bf Example 2.2:} Let $R = A[X_1, X_2, \ldots, X_t]$ be the ring of polynomial in $t$
 indeterminates $X_1, X_2, \ldots, X_t$ with coefficients in $A$ ($\dim A = d > 0$). Then
$R =\bigoplus_{m \geq 0}R_m$ 
is a finitely generated standard graded algebra over $A$, where $R_m$
is the set of all homogeneous polynomials of degree $m$ and the zero polynomial. It is well known that $X_1, X_2, \ldots, X_t$ is a regular sequence of $R$. Let $I$ be an ideal of $A$. It is easy to see that  $X_1R_{m-1} \cap IR_m $ and $IX_1R_{m-1}$ are both  the set of all homogeneous polynomials of degree $m$ with coefficients in $\frak{I}$ and divided by $X_1$. Hence $X_1R_{m-1} \cap IR_m = IX_1R_{m-1}$ for any ideal $I$ of $A$. Using the results just obtained and the fact that 
$$R/(X_1, \ldots, X_i)R = A[X_{i + 1}, \ldots, X_t]$$
for all $i < t$, we immediately show that
 $X_1, X_2, \ldots, X_t$ be a weak-(FC)-sequence of $R$ with respect to $I$ for any ideal $I$ of $A$. 
 \vskip 0.2cm
 Now, we give some comments on weak-(FC)-sequences of a finitely generated standard multi-graded algebra over $A$ by the following remark.

\noindent {\bf Remark 2.3.} 
\begin{itemize}
\item[(i)]  By Artin-Rees lemma, there exists integer $u_1, u_2, \ldots, u_s$  such that 
\begin{align*}
(0 : S_{++}^\infty) \bigcap S_{(n_1, \ldots, n_s)}& = S_{(n_1- u_1, \ldots, n_s - u_s)}((0 : S_{++}^\infty) \bigcap S_{(u_1, \ldots, u_s)})\\
& \subseteq S_{(n_1- u_1, \ldots, n_s - u_s)}(0 : S_{++}^\infty) 
\end{align*}
for all $n_1\geq  u_1, \ldots, n_s \geq  u_s$. Since 
$S_{(n_1- u_1, \ldots, n_s - u_s)}(0 : S_{++}^\infty) = 0$
for all large enough $n_1, \ldots, n_s$, it follows that
$(0 : S_{++}^\infty)_{(n_1, \ldots, n_s)} = (0 : S_{++}^\infty) \bigcap S_{(n_1, \ldots, n_s)} = 0$
for all large enough $n_1, \ldots, n_s$. 
\item[(ii)]  Let $x \in S$ be a homogeneous element. If $0 : x \subseteq  0: S_{++}^\infty $ then, by (i), $$(0 : x)_{(n_1, \ldots, n_s)}\subseteq (0 : S_{++}^\infty)_{(n_1, \ldots, n_s)} = 0$$ for all large $n_1, \ldots, n_s$. Thus $x$ is a filter-regular element of $S$.  Conversely, suppose that $x$ is a filter-regular element of $S$.  We have   
$$S_{(n_1, \ldots, n_s)}(0 : x)_{(v_1, \ldots, v_s)} \subseteq (0 : x)_{(n_1 + v_1, \ldots, n_s + v_s)} = 0  $$
for all large $n_1, \ldots, n_s $ and all $v_1, \ldots, v_s$. It implies that 
$$(0 : x)_{(v_1, \ldots, v_s)} \subseteq  (0 : S_{++}^n) \subseteq (0 : S_{++}^\infty )$$
 for all large $n$ and all $v_1, \ldots, v_s$.  Hence 
$(0 : x) \subseteq (0 : S_{++}^\infty )$.  Therefore $x$ is a filter-regular element of $S$ if and only if $0 : x \subseteq  0: S_{++}^\infty $. 
\item[(iii)] Suppose that $x \in S_i$ is a filter-regular element of $S$. Consider
$$\lambda_x : S_{(n_1, \ldots, n_i, \ldots, n_s)} \longrightarrow xS_{(n_1, \ldots, n_i - 1, \ldots, n_s)}, y \mapsto xy.$$
It is clear that $\lambda_x $ is surjective and 
$\text{ker}\lambda_x = (0 : x) \cap S_{(n_1, \ldots, n_s)} = 0$ 
for all large $n_1, \ldots, n_s$. Therefore, 
$S_{(n_1, \ldots, n_i, \ldots, n_s)} \cong  xS_{(n_1, \ldots, n_i - 1, \ldots, n_s)}$. This follows that
$$IS_{(n_1, \ldots, n_i, \ldots, n_s)} \cong  xIS_{(n_1, \ldots, n_i - 1, \ldots, n_s)}$$
for all large $n_1, \ldots, n_s$ and for any ideal $I$ of $A$. 
\item[(iv)] If $S_{++}$ is non-nilpotent then $S_{(n, \ldots, n)} \not= 0$ for all $n$. Hence, by Nakayama's Lemma, 
$(D_\frak{m}(S))_n = \frac{\frak{m}^nS_{(n, \ldots, n)}}{\frak{m}^{n + 1}S_{(n, \ldots, n)}} \not = 0$
for all $n$. It implies that $\dim D_\frak{m}(S) \geq 1$.
\end{itemize}

\vskip 0.2cm
The following lemma will play a crucial role for showing the existence of weak-(FC)-sequence. 
\vskip 0.2cm
\noindent {\bf Lemma 2.4 (Generalized Rees's Lemma).} {\it Let $(A, \mathfrak{m})$ be a Noetherian local ring with maximal ideal $\mathfrak{m}$, infinite residue $k = A/\mathfrak{m}$. 
Let $S= \bigoplus_{n_1, \ldots, n_s \geq 0}S_{(n_1, \ldots,  n_s)}$ 
be a finitely generated standard s-graded algebra over $A$; $I$ be an ideal of $A$. 
Let $\Sigma$ be a finite collection of prime ideals of $S$ not containing $S_{(1, \ldots, 1)}$.
Then for each $i = 1, \ldots , s$, there exists an element $x_i \in S_i \backslash \mathfrak{m}S_i$,
$x_i$ not contained in any prime ideal in $\Sigma$, and a positive integer $k_i$ such that
$$x_iS_{(n_1, ..., n_i - 1, ..., n_s)}\cap I^{n_0}S_{(n_1, ..., n_s)} 
= x_iI^{n_0}S_{(n_1, \ldots, n_i - 1, \ldots, n_s)}$$ 
for all $n_i > k_i$ and all non-negative integers $n_0, n_1, \ldots, n_{i-1}, n_{i+1}, \ldots, n_s$.}
\vskip 0.2cm
 \noindent{\bfseries\itshape\textit Proof.} In the ring $S[t, t^{-1}]$ ($t$ is an indeterminate), set
$$S^* = \bigoplus_{n_0 \in \mathbb{Z}}I^{n_0}St^{n_0} = 
\bigoplus_{n_0 \in \mathbb{Z}; n_1, \ldots, n_s \geq 0}I^{n_0}S_{(n_1, \ldots, n_s)}t^{n_0} $$
where $I^n = A$ for $n \leq 0$. Then $S^*$ is a Noetherian $(s + 1)$-graded ring. 
From $u=t^{-1}$ is  non-zero-divisor in $S^*$, by the Corollary of [Lemma 2.7, Re], 
the set of prime associated with $u^nS^*$ is independent on $n > 0$ and so is finite. 
We divide this set into two subsets: $\mathfrak{S}_1$ consisting of those containing $S_i$ and 
$\mathfrak{S}_2$ those that do not
(where $S_i = S_{(0, \ldots, \underset{(i)}{1}, \ldots, 0)} = S^*_{(0, 0, \ldots, \underset{(i + 1)}{1}, \ldots, 0)}$).\\
From $S_i/\mathfrak{m}S_i$ is a vector space over the infinite field $k$ and the sets 
$\Sigma, \mathfrak{S_2}$ are both finite, we can choose 
$x_i \in S_i \backslash \mathfrak{m}S_i$ such that $x_i$ is not contained in any prime ideal belonging to $\Sigma \cup  \mathfrak{S_2}$.
Set $$M_n = \frac{(u^nS^* : x_i) \cap S^*}{u^nS^*}.$$
Then $M_n$ is a $S^*$-module for any $n > 0$. We need must show that there exists a sufficiently large integer $N > 0$ such that $S_i^NM_n = 0$. Note that if $P \in  \text{Ass}_{S^*}M_n$ then $P \in \text{Ass}_{S^*}S^*/u^nS^* = \mathfrak{S}_1\cup \mathfrak{S}_2$, and there exists $z \in u^nS^* : x_i$ such that $P= u^nS^* : z$. Since $x_iz \in u^nS^*$, $x_i \in P$. So $P \in \mathfrak{S}_1$. Hence $S_i \subset  P$. It follows that 
$S_i \subset  \bigcap_{P\in \text{Ass}_{S^*}M_n} P$. Therefore
$$S_i \subset  \sqrt{\text{Ann}_{S^*}M_n}.$$
Since $S_i$ is finitely generated, there exists a sufficiently large integer $N > 0$ 
(how large depending on $n$) such that $S_i^NM_n = 0$. Hence,  for all large $n_i > N$, 
any element of $M_n$ of degree $(n_0, n_1, \ldots, n_s)$ is zero. This means that, for each $n > 0$, we have
\begin{equation}\label{E: pt1}
(u^nI^{n_0}S_{(n_1, ..., n_s)}t^{n_0} : x_i)\bigcap S^*  = u^nI^{n_0}S_{(n_1, ...,n_i - 1, ..., n_s)}t^{n_0} 
\end{equation}
for all large $n_i$ and all non-negative integers $n_0, n_1, \ldots, n_{i - 1}, n_{i+1}, \ldots, n_s$.\\
Let $\mathfrak{b}$ denote the ideal of $S^*$ consisting of all finite sums $\sum c_{n_0}t^{n_0}$ with
$$c_{n_0} \in x_iS_{(n_1, ..., n_i - 1, ..., n_s)}\cap I^{n_0}S_{(n_1, ..., n_s)}.$$
Then $\mathfrak{b}$  has a finite generating set of the form 
$x_ib_i t^{n_0}$ with $b_i \in S_{(n_1, \ldots, n_i - 1, \ldots, n_s)}$.  
Note that if $0 \not= a \in I^mS$ and $m \geq n_0$ then $at^{n_0} \in S^*$. 
Specially, if $n_0 < 0$  then $at^{n_0} \in S^*$ for all $a \in S$.  Hence since the above generating set of $\mathfrak{b}$  is finite,  it follows that there exists an integer $q$ such that
$u^qb_it^{n_0} = b_it^{n_0 - q} \in S^*$ for all element of this generating set ($q$ is chosen such that  $n_0 - q < 0$ for all $n_0$). Therefore $\mathfrak{b} \subseteq x_iS^* : u^q$.\\
Now, suppose that $z \in x_iS_{(n_1, ..., n_i - 1, ..., n_s)}\cap I^{n_0}S_{(n_1, ..., n_s)}$.
This means $zt^{n_0} \in \mathfrak{b}$. Because  $\mathfrak{b} \subseteq x_iS^* : u^q$,  $u^qzt^{n_0}  = x_iw$, where $w \in S^*$.
Since $z \in I^{n_0}S_{(n_1, ..., n_s)}$, it follows that $x_iw = u^qzt^{n_0} \in u^qI^{n_0}S_{(n_1, ..., n_s)}t^{n_0}$.
Hence, by (\ref{E: pt1}), we can find $k_i$ such that 
 $$w \in (u^qI^{n_0}S_{(n_1, ..., n_s)}t^{n_0} : x_i)\bigcap S^*
 = u^qI^{n_0}S_{(n_1, ...,n_i - 1, ..., n_s)}t^{n_0}$$
for all $n_i > k_i $. Thus
$u^qzt^{n_0} = x_iw \in x_iu^qI^{n_0}S_{(n_1, ...,n_i - 1, ..., n_s)}t^{n_0}.$
Since $u$ and $t$ are non-zero-divisors in $S^*$,
 $z \in x_iI^{n_0}S_{(n_1, ...,n_i - 1, ..., n_s)}.$
Hence if $n_i > k_i$,
$$x_iS_{(n_1, ..., n_i - 1, ..., n_s)}\cap I^{n_0}S_{(n_1, ..., n_s)} 
\subset  x_iI^{n_0}S_{(n_1, ...,n_i - 1, ..., n_s)}.$$
Consequently,
$x_iS_{(n_1, ..., n_i - 1, ..., n_s)}\cap I^{n_0}S_{(n_1, ..., n_s)} =  x_iI^{n_0}S_{(n_1, ...,n_i - 1, ..., n_s)}.  \ \ \blacksquare$ 
 \vskip 0.2cm
The following proposition will show the existence of weak-(FC) sequence.
 \vskip 0.2cm
\noindent {\bf Proposition 2.5.} {\it  Suppose that $S_{++}$ is non-nilpotent.
Then for any $1\leq i \leq s$, there exists a weak-(FC)-element $x \in S_i$ of $S$ with respect to $I$.}
\vskip 0.2cm
  \noindent{\bfseries\itshape\textit Proof.} Since $S_{++}$ is non-nilpotent,
 $ S/0 : S_{++}^\infty \not= 0$.  Set 
$$\Sigma = \text{Ass}_S(S/0: S_{++}^\infty) = \{P \in \text{Ass}S \mid P \nsupseteq  S_{(1, \ldots, 1)}\}.$$
 Then $\Sigma$ is finite. By Lemma 2.4, for each $i = 1, \ldots, s$, there exists 
$x \in S_i \setminus \mathfrak{m}S_i$ such that $x \notin P$ for all $P \in \Sigma$ and 
$$xS_{(n_1, ..., n_i - 1, ..., n_s)}\cap I^{n_0}S_{(n_1, \ldots, n_i - 1, \ldots, n_s)} 
= xI^{n_0}S_{(n_1, ..., n_s)}.$$ 
Thus $x$ satisfies the condition (FC$_1$).
Since $x \notin P$ for all $P \in \Sigma$,  $0 : x \subset 0 : S_{++}^\infty $. Hence by Remark 2.3(ii), $x$ satisfies the condition (FC$_2$). $\blacksquare$
 \vskip 0.2cm
\noindent {\bf 3. Mixed multiplicities of multi-graded algebras} 
 
 This section first determines mixed multiplicities of multi-graded algebras defined over a Noetherian local ring, next answers to the question when these mixed multiplicities are positive and characterizes them in term of Hilbert-Samuel multiplicities.

Let $S= \bigoplus_{n_1, \ldots, n_s \geq 0}S_{(n_1, \ldots, n_s)}$ 
be a finitely generated standard $s$-graded algebra over a Noetherian local ring $A$
such that $S_{++}$ is non-nilpotent and $J$ be an $\mathfrak{m}$-primary ideal of $A$. Since 
$$\bigoplus_{n_0, n_1, \ldots, n_s \geq 0}\frac{J^{n_0}S_{(n_1, \ldots, n_s)}}{J^{n_0 + 1}S_{(n_1, \ldots, n_s)}}$$
is a finitely generated standard $s$-graded algebra over Artinian local ring $A/J$, by [HHRT, Theorem 4.1], 
 $$\ell_A(\frac{J^{n_0}S_{(n_1, \ldots, n_s)}}{J^{n_0 + 1}S_{(n_1, \ldots, n_s)}})$$
 is a polynomial for all large $n_0, n_1, \ldots, n_s$. Denote by $P(n_0, n_1, \ldots, n_s)$ this 
polynomial.
Set $$D_J(S) = \bigoplus_{n \geq  0} \frac{J^nS_{(n, ..., n)}}{J^{n + 1}S_{(n, ..., n)}}$$
and $\ell = \dim{D_\frak{m}(S)}$. By Remark 2.3(iv), $\ell \geq 1$. Note that $\dim D_J(S) = \dim{D_\frak{m}(S)}$ for all $\frak{m}$-primary ideal $J$ of $A$ and 
$\deg P(n_0, n_1, \ldots, n_s) = \deg P(n, n, \ldots, n)$.  Since
$$P(n, n, \ldots, n) = \ell_A(\frac{J^{n}S_{(n, \ldots, n)}}{J^{n + 1}S_{(n, \ldots, n)}}) = \ell_A(D_J(S)_n)$$ 
for all large $n$, it follows that 
$\deg P(n, n, \ldots, n) = \dim D_J(S) - 1 = \ell - 1.$
Hence $\deg P(n_0, n_1, \ldots, n_s) = \ell - 1$. 

If we write the term of total degree $\ell - 1$ of $P$ in the form
$$\sum_{k_0 + k_1 + \cdots + k_s = \ell - 1} e(J, k_0, k_1, \ldots, k_s, S) 
\frac{n_0^{k_0}n_1^{k_1}\cdots n_s^{k_s}}{n_0!n_1! \cdots n_s!}$$
then $e(J, k_0, k_1, \ldots, k_s, S)$ are non-negative integers not all zero and called the 
{\it mixed multiplicity of $S$ of type $(k_0, k_1, \ldots, k_s)$ with respect to $J$}.
 
From now on, the notation $e_A(J, M)$ will mean the Hilbert-Samuel multiplicity of $A$-module $M$ with respect to ideal $\frak{m}$-primary $J$ of $A$.  We shall begin the section with the following lemma.
  \vskip 0.2cm
 \noindent {\bf Lemma 3.1.} {\it Let $S$ be a finitely generated standard $s$-graded algebra over a Noetherian local ring $A$ such that $S_{++}$ is non-nilpotent and $J$ be an 
$\mathfrak{m}$-primary ideal of $A$. Set $\ell = \dim D_\frak{m}(S)$. Then
$e(J, k_0, 0, \ldots, 0, S) \not= 0 $ if and only if
$\dim A/(0 : S_{(1, \ldots, 1)}^\infty)_A = \ell.$
 In this case, $e(J, k_0, 0, \ldots, 0, S) = e_A(J, S_{(n, \ldots, n)})$ for all large $n$.}
\vskip 0.2cm
  \noindent{\bfseries\itshape\textit Proof.}
  Denote by $P(n_0, n_1, \ldots, n_s)$ the polynomial of
$$\ell_A(\frac{J^{n_0}S_{(n_1, \ldots, n_s)}}{J^{n_0 + 1}S_{(n_1, \ldots, n_s)}}).$$
Then $P$ is a polynomial of total degree $\ell - 1$. By taking 
$n_1 = n_2 = \cdots = n_s = u$, where $u$ is a sufficiently large integer, we get
$$e(J, k_0, 0, \ldots, 0, S) = \lim\limits_{n_0\to\infty }
\frac{(\ell - 1)!P(n_0, u, \ldots, u) }{n_0^{\ell - 1}}.$$
Since $P(n_0, u, \ldots, u) = \ell_A(\frac{J^{n_0}S_{(u, \ldots, u)}}{J^{n_0 + 1}S_{(u, \ldots, u)}})$,
  it follows that
   
    \centerline {$\deg P(n_0, u, \ldots, u) = \dim_AS_{(u, \ldots, u)} - 1$}  
\noindent and 
$e(J, k_0, 0, \ldots, 0, S) \not= 0$
if and only if  $$\deg P(n_0, u, \ldots, u) = \dim_AS_{(u, \ldots, u)} - 1 = \ell - 1.$$
Since $A$ is Noetherian, 
$(0 : S_{(1, \ldots, 1)}^\infty)_A = (0: S_{(1, \ldots, 1)}^n)_A = (0 : S_{(n, \ldots, n)})_A$
\noindent for all large $n$. Hence if $u$ is chosen sufficiently large, we have
 
 \centerline {$\dim_AS_{(u, \ldots, u)} = \dim A/(0 : S_{(u, \ldots, u)})_A 
= \dim A/(0 : S_{(1, \ldots, 1)}^\infty)_A .$}
 \noindent
 Therefore
 $e(J, k_0, 0, \ldots, 0, S)  \not= 0$  if and only if   $\dim A/(0 : S_{(1, \ldots, 1)}^\infty)_A  = \ell .$ 
Finally, if  $\dim A/(0 : S_{(1, \ldots, 1)}^\infty)_A  = \ell $ then
$\dim_AS_{(n, \ldots, n)} - 1= \ell - 1$ for all large $n$ and hence
\begin{align*}
e_A(J, S_{(n, \ldots, n)}) &= 
\lim\limits_{n_0\to\infty }
\frac{(\ell  - 1)!\ell_A(\frac{J^{n_0}S_{(n, \ldots, n)}}{J^{n_0 + 1}S_{(n, \ldots, n)}})}{n_0^{\ell  - 1}}\\
&= \lim\limits_{n_0\to\infty }
\frac{(\ell -  1)!P(n_0, n, \ldots, n) }{n_0^{\ell - 1}} 
= e(J, k_0, 0, \ldots, 0, S) 
\end{align*}
for all large integer $n$. $\blacksquare$
 \vskip 0.3cm
 \noindent {\bf Proposition 3.2.} {\it Let $S$ be a finitely generated standard $s$-graded algebra over 
a Noetherian local ring $A$ such that $S_{++}$ is non-nilpotent and $J$ be an  
$\mathfrak{m}$-primary ideal of $A$. Set $\ell = \dim D_\frak{m}(S)$.
Assume that  $e(J, k_0, k_1, \ldots, k_s, S) \not= 0$, where $k_0, k_1, \ldots, k_s$ are non-negative  integers such that $k_0 + k_1 + \cdots + k_s = \ell - 1$. Then 
 \begin{enumerate}
 \item[\rm (i)] If  $k_i > 0$ and $x \in S_i$ is a weak-(FC)-element of $S$ with respect to $J$ then
$$e(J, k_0, k_1, \ldots, k_s, S) = e(J, k_0, k_1, \ldots, k_i - 1, \ldots, k_s, S/xS),$$
and $\dim D_\frak{m}(S/xS) = \ell - 1$.
 \item[\rm (ii)] There exists a weak-(FC)-sequence of $t = k_1 + \cdots + k_s$ elements of $S$
 in $\bigcup_{i = 1}^sS_i$  with respect to $J$
 consisting of $k_1$ elements of $S_1$, $\ldots$, $k_s$ elements of $S_s$.
\end{enumerate} } 
\vskip 0.2cm
  \noindent{\bfseries\itshape\textit Proof.} 
 The proof of (i): Denote by $P(n_0, n_1, \ldots , n_s)$ the polynomial of 
$$\ell_A(\frac{J^{n_0}S_{(n_1, ..., n_s)}}{J^{n_0 + 1}S_{(n_1, ..., n_s)}}).$$
Then $\deg P = \ell - 1$. Since $x$ satisfies the condition (FC$_1$), for all large $n_0, n_1, \ldots, n_s$, we have 
\begin{align*}
&\ell_A(\frac{J^{n_0}(S/xS)_{(n_1, \ldots, n_s)}}{J^{n_0 + 1}(S/xS)_{(n_1, \ldots, n_s)}})
 = \ell_A(\frac{J^{n_0}(S_{(n_1, \ldots, n_s)}\big/xS_{(n_1, \ldots, n_i - 1, \ldots, n_s)})}
                     {J^{n_0 + 1}(S_{(n_1, \ldots, n_s)}\big/xS_{(n_1, \ldots, n_i - 1, \ldots, n_s)})}) \\\
& = \ell_A(\frac{J^{n_0}S_{(n_1, \ldots, n_s)} + xS_{(n_1, \ldots, n_i - 1, \ldots, n_s)}}
                     {J^{n_0 + 1}S_{(n_1, \ldots, n_s)} + xS_{(n_1, \ldots, n_i - 1, \ldots, n_s)}}) \\\
& = \ell_A(\frac{J^{n_0}S_{(n_1, \ldots, n_s)}}{(J^{n_0 + 1}S_{(n_1, \ldots, n_s)} + 
                                xS_{(n_1, \ldots, n_i - 1, \ldots, n_s)})\cap J^{n_0}S_{(n_1, \ldots, n_s)}}) \\\                     
& = \ell_A(\frac{J^{n_0}S_{(n_1, \ldots, n_s)}}{J^{n_0 + 1}S_{(n_1, \ldots, n_s)} + 
                                xS_{(n_1, \ldots, n_i - 1, \ldots, n_s)}\cap J^{n_0}S_{(n_1, \ldots, n_s)}}) \\\
& = \ell_A(\frac{J^{n_0}S_{(n_1, \ldots, n_s)}}
                {J^{n_0 + 1}S_{(n_1, \ldots, n_s)} +  x J^{n_0}S_{(n_1, \ldots, n_i - 1, \ldots, n_s)}}) \\\    
& = \ell_A(\frac{J^{n_0}S_{(n_1, \ldots, n_s)}}{J^{n_0 + 1}S_{(n_1, \ldots, n_s)}})
   - \ell_A(\frac{J^{n_0 + 1}S_{(n_1, \ldots, n_s)} +  x J^{n_0}S_{(n_1, \ldots, n_i - 1, \ldots, n_s)}}
                   {J^{n_0 + 1}S_{(n_1, \ldots, n_s)}}) \\\
& = \ell_A(\frac{J^{n_0}S_{(n_1, \ldots, n_s)}}{J^{n_0 + 1}S_{(n_1, \ldots, n_s)}})
      - \ell_A(\frac{x J^{n_0}S_{(n_1, \ldots, n_i - 1, \ldots, n_s)}}
            {x J^{n_0}S_{(n_1, \ldots, n_i - 1, \ldots, n_s)}\cap J^{n_0 + 1}S_{(n_1, \ldots, n_s)}}) \\\
&  = \ell_A(\frac{J^{n_0}S_{(n_1, \ldots, n_s)}}{J^{n_0 + 1}S_{(n_1, \ldots, n_s)}})
      - \ell_A(\frac{x J^{n_0}S_{(n_1, \ldots, n_i - 1, \ldots, n_s)}}
            {x J^{n_0 + 1}S_{(n_1, \ldots, n_i - 1, \ldots, n_s)}}). 
\end{align*}\\
Since $x$ is a filter-regular element of $S$, it follows by Remark 2.3(iii) that
$$J^{n_0}S_{(n_1, \ldots, n_i, \ldots, n_s)} \cong  xJ^{n_0}S_{(n_1, \ldots, n_i - 1, \ldots, n_s)}$$
for all $n_0$ and all large $n_1, \ldots, n_s$.  Thus we have an isomorphism of $A$-modules
$$\frac{x J^{n_0}S_{(n_1, \ldots, n_i - 1, \ldots, n_s)}}
            {x J^{n_0 + 1}S_{(n_1, \ldots, n_i - 1, \ldots, n_s)}} \cong  
\frac{J^{n_0}S_{(n_1, \ldots, n_i - 1, \ldots, n_s)}}
            {J^{n_0 + 1}S_{(n_1, \ldots, n_i - 1, \ldots, n_s)}}$$
for all large $n_0, n_1, \ldots, n_s$. So
$$\ell_A(\frac{x J^{n_0}S_{(n_1, \ldots, n_i - 1, \ldots, n_s)}}
            {x J^{n_0 + 1}S_{(n_1, \ldots, n_i - 1, \ldots, n_s)}}) = 
\ell_A(\frac{J^{n_0}S_{(n_1, \ldots, n_i - 1, \ldots, n_s)}}
            {J^{n_0 + 1}S_{(n_1, \ldots, n_i - 1, \ldots, n_s)}}) .$$
Hence
$$\ell_A(\frac{J^{n_0}(S/xS)_{(n_1, \ldots, n_s)}}{J^{n_0 + 1}(S/xS)_{(n_1, \ldots, n_s)}}) = 
     \ell_A(\frac{J^{n_0}S_{(n_1, \ldots, n_s)}}{J^{n_0 + 1}S_{(n_1, \ldots, n_s)}})
    - \ell_A(\frac{J^{n_0}S_{(n_1, \ldots, n_i - 1, \ldots, n_s)}}
                       {J^{n_0 + 1}S_{(n_1, \ldots, n_i - 1, \ldots, n_s)}})$$
for all large $n_0, n_1, \ldots, n_s$.
Denote by $Q(n_0, n_1, \ldots, n_s)$ the polynomial of
$$\ell_A(\frac{J^{n_0}(S/xS)_{(n_1, \ldots, n_s)}}{J^{n_0 + 1}(S/xS)_{(n_1, \ldots, n_s)}}).$$
From the above fact, we get

\centerline {$Q(n_0, n_1, \ldots, n_s) = P(n_0, n_1, \ldots, n_i, \ldots, n_s) - P(n_0, n_1, \ldots, n_i - 1, \ldots, n_s).$}
 \noindent Since $e(J, k_0, k_1, \ldots,  k_s, S) \not= 0 $ and $k_i > 0$, it implies that $\deg Q =  \deg P - 1$ and $$e(J, k_0, k_1, \ldots, k_i , \ldots, k_s, S) = e(J, k_0, k_1, \ldots, k_i - 1, \ldots, k_s, {S/xS}).$$
 Note that  $\deg Q = \dim D_\frak{m}(S/xS) - 1$. Hence
 $$\dim D_\frak{m}(S/xS) = \deg Q + 1 = \deg P = \ell - 1.$$
The proof of (ii): 
 The proof is by induction on $t = k_1 + \cdots + k_s$. For $t  = 0$, the result is trivial. 
Assume that $t > 0$. Since $k_1 + \cdots + k_s = t > 0$, there exists $k_j > 0$. 
 Since $S_{++}$ is non-nilpotent, by Proposition 2.5, there exists a weak-(FC) element $x_1 \in S_j$ of $S$ with respect to $J$. 
By (i), 
$$e(J, k_0, k_1, \ldots, k_i - 1, \ldots, k_s, S/x_1S) = e(J, k_0, k_1, \ldots, k_s, S) \not= 0.$$
This follows that $$\frac{J^{n_0}(S/x_1S)_{(n_1, \ldots, n_s)}}{J^{n_0 + 1}(S/x_1S)_{(n_1, \ldots, n_s)}}$$ and so $(S/x_1S)_{(n_1, \ldots, n_s)}\ne 0$  for all large $n_1, \ldots, n_s$. Hence $(S/x_1S)_{++}$ is non-nilpotent. 
Since $k_1 + \cdots + k_j - 1 + \cdots + k_s = t - 1$, by the inductive assumption, there exists 
$t - 1$ elements $x_2, \ldots, x_t$ consisting of $k_1$ elements of $S_1$, ...,  $k_j - 1$ elements 
of $S_j$, ..., $k_s$ elements of $S_s$ such that $\bar{x}_2, \ldots, \bar{x}_t$ is 
a weak-(FC)-sequence  of $S/x_1S$ with respect to $J$ ($\bar{x}_i$ is initial form of $x_i$ in 
$S/x_1S$, $i = 2, \ldots, t$). Hence $x_1, \ldots, x_t$ is a weak-(FC)-sequence of $S$ with respect
 to $J$ consisting of $k_1$ elements of $S_1$, ..., $k_s$ elements of $S_s$. $\blacksquare$
  \vskip 0.2cm
The following theorem will give the criteria for the positivity of mixed multiplicities and characterize them in term of Hilbert-Samuel multiplicity.
 \vskip 0.2cm
\noindent {\bf Theorem 3.3.} {\it Let $S$ be a finitely generated standard $s$-graded algebra over a Noetherian local ring $A$ such that $S_{++}$ is non-nilpotent. Let $J$ be an $\frak{m}$-primary ideal of $A$. Set $\ell = \dim D_\frak{m}(S)$. Then the following statements hold.
\begin{enumerate}
\item[\rm (i)] $e(J, k_0, k_1, \ldots, k_s, S) \not= 0$ if and only if there exists a weak-(FC)-sequence 
$x_1,  \ldots, x_t$ $(t = k_1 + \cdots + k_s)$ of $S$ with respect to $J$ 
consisting of $k_1$ elements of $S_1$, ..., $k_s$ elements of $S_s$ and 
$$\dim D_\frak{m}(S/(x_1, \ldots, x_t)S) = \dim A/((x_1, \ldots, x_t)S: {S}_{(1, \ldots, 1)}^\infty )_A= \ell - t.$$
\item[\rm (ii)] Suppose that $e(J, k_0, k_1, \ldots, k_s, S) \not= 0$ and $x_1,  \ldots, x_t$ 
$(t = k_1 + \cdots + k_s)$ is a weak-(FC)-sequence of $S$ with respect to $J$ consisting of $k_1$ elements of $S_1$, ..., $k_s$ elements of $S_s$. Set $\bar{S} = S/(x_1, \ldots, x_t)S$. Then
$$e(J, k_0, k_1, \ldots, k_s, S) = e_A(J, \bar{S}_{(n, \ldots, n)})$$
for all large $n$.
\end{enumerate}}
\vskip 0.2cm
 \noindent{\bfseries\itshape\textit Proof.}
The proof of (i):
First, we prove the necessary condition. By Proposition 3.2(ii), there exists a weak-(FC)-sequence 
$x_1, \ldots, x_t$ of $S$ with respect to $J$ consisting of $k_1$ elements of $S_1$, ..., $k_s$ 
elements of $S_s$. Set $\bar{S} = S/(x_1, \ldots, x_t)S$. Applying Proposition 3.2(i) by induction on $t$, we get $\dim D_\frak{m}(\bar{S}) = \ell - t$ and
$$0 \not= e(J, k_0, k_1, \ldots, k_s, S) = e(J, k_0, 0, \ldots, 0, \bar{S}).$$
Hence by Lemma 3.1, 
$\dim A/(0 : \bar{S}_{(1, \ldots, 1)}^\infty)_A = \ell - t.$
Since
 $$\dim A/(0 : \bar{S}_{(1, \ldots, 1)}^\infty)_A = \dim A/((x_1, \ldots, x_t)S: {S}_{(1, \ldots, 1)}^\infty )_A,$$
 it follows that
$\dim D_\frak{m}(S/(x_1, \ldots, x_t)S) = \dim A/((x_1, \ldots, x_t)S: {S}_{(1, \ldots, 1)}^\infty )_A= \ell - t.$
Now, we prove the sufficiently condition. Without loss of general, we may assume that $x_1 \in S_i$.
Denote by $P(n_0, n_1, \ldots, n_s)$ and $Q(n_0, n_1, \ldots, n_s)$ the polynomials of
$$\ell_A(\frac{J^{n_0}S_{(n_1, \ldots, n_s)}}{J^{n_0 + 1}S_{(n_1, \ldots, n_s)}}) \ 
\text{ and }
 \ \ell_A(\frac{J^{n_0}(S/x_1S)_{(n_1, \ldots, n_s)}}{J^{n_0 + 1}(S/x_1S)_{(n_1, \ldots, n_s)}}),$$
 respectively. Then by the proof of Proposition 3.2(i) we have
$$Q(n_0, n_1, \ldots, n_s) = P(n_0, n_1, \ldots, n_i, \ldots, n_s) - P(n_0, n_1, \ldots, n_i - 1, \ldots, n_s).$$
This implies that $\deg Q \leq \deg P - 1$. Recall that $\deg Q = \dim D_\frak{m}(S/x_1S) - 1$ and $\deg P = \dim D_\frak{m}(S) - 1$. So $\dim D_\frak{m}(S/x_1S) \leq  \dim D_\frak{m}(S) - 1$. Similarly, we have
\begin{align*}
\ell - t =& \dim D_\frak{m}(S/(x_1, \ldots, x_t)S) \leq \dim D_\frak{m}(S/(x_1, \ldots, x_{t- 1})S) - 1\\
& \leq \cdots \leq \dim D_\frak{m}(S/x_1S) - (t - 1) \leq \dim D_\frak{m}(S) - t = \ell - t . 
\end{align*}
This fact follows $\dim D_\frak{m}(S/x_1S) = \dim D_\frak{m}(S) - 1$.  Thus $\deg Q = \deg P - 1$. Hence
$$e(J, k_0, k_1, \ldots, k_s, S) = e(J, k_0, k_1, \ldots, k_i - 1, \ldots, k_s, S/x_1S).$$
By induction we have
$e(J, k_0, k_1, \ldots, k_s, S) = e(J, k_0, 0, \ldots, 0, \bar{S}).$
Since $$\dim A/(0 : \bar{S}_{(1, \ldots, 1)}^\infty)_A = \dim A/((x_1, \ldots, x_t)S: {S}_{(1, \ldots, 1)}^\infty )_A  = \ell - t = \dim D_\frak{m}(\bar{S}),$$ it follows,  by Lemma 3.1, that
$e(J, k_0, 0, \ldots, 0, \bar{S}) \not= 0$.
Hence $$e(J, k_0, k_1, \ldots, k_s, S) \not= 0.$$

The proof of (ii): Applying Proposition 3.2(i), by induction on $t$, we obtain
$$0 \not= e(J, k_0, k_1, \ldots, k_s, S) = e(J, k_0, 0, \ldots, 0, \bar{S}).$$
On the other hand, by Lemma 3.1, 
$e(J, k_0, 0, \ldots, 0, \bar{S}) = e_A(J, \bar{S}_{(n, \ldots, n)})$
for all large integer $n$. Hence $e(J, k_0, k_1, \ldots, k_s, S) = e_A(J, \bar{S}_{(n, \ldots, n)})$ for all large $n$.
$\blacksquare$
\vskip 0.2 cm
\noindent{\bf Remark 3.4.} From the proof of Theorem 3.3 we get  some comments as following. 
 \begin{itemize}
 \item[\rm (i)] Assume that $x_1, \ldots, x_t$ is a weak-(FC)-sequence in $\bigcup_{i = 1}^sS_i$ of $S$ with respect to $J$. If  $\dim D_\frak{m}(S/(x_1, \ldots, x_t)S) = \dim D_\frak{m}(S) - t$ then $\dim D_\frak{m}(S/(x_1, \ldots, x_i)S) = \dim D_\frak{m}(S) - i$ for all $ 1 \leq i \leq t$. 
\item[\rm (ii)] If  $k_i > 0$ and $x \in S_i$  is a weak-(FC)-sequence of $S$ with respect to $J$ such that $\dim D_\frak{m}(S/xS) = \dim D_\frak{m}(S) - 1$ then 
 $$e(J, k_0, k_1, \ldots, k_i , \ldots, k_s, S) = e(J, k_0, k_1, \ldots, k_i - 1, \ldots, k_s, S/xS).$$
\item[\rm (iii)] If $e(J, k_0, k_1, \ldots, k_s, S) \not= 0$ then for every weak-(FC)-sequence $x_1, \ldots, x_t$ ($t = k_1 + \cdots + k_s$) of $S$  with respect  to $J$ consisting of $k_1$ elements of $S_1$, ..., $k_s$ elements of $S_s$ we always have
  $$\dim D_\frak{m}(S/(x_1, \ldots, x_t)S) = \dim A/((x_1, \ldots, x_t)S: {S}_{(1, \ldots, 1)}^\infty )_A= \ell - t.$$
 \item[\rm (iv)] Suppose that $x_1, \ldots, x_t$ is a weak-(FC)-sequence in $\bigcup_{i = 1}^sS_i$ of $S$ with respect to $J$. Then $\dim D_\frak{m}(S/x_1S) \leq  \dim D_\frak{m}(S) - 1$. By induction we have 
$$\dim D_\frak{m}(S/(x_1, \ldots, x_t)S) \leq  \dim D_\frak{m}(S) - t = \ell - t.$$ 
If $\ell = t$ then $\dim D_\frak{m}(S/(x_1, \ldots, x_t)S) \leq 0$. Hence $(S/(x_1, \ldots, x_t)S)_{++}$ is nilpotent by Remark 2.3(iv). So $x_1, \ldots, x_t$ is a maximal weak-(FC)-sequence. This fact follows that the length of every weak-(FC)-sequence in $\bigcup_{i = 1}^sS_i$ of $S$ with respect to $J$ is not greater than $\ell$.
  \end{itemize}

\vskip 0.2cm
\noindent {\bf Example 3.5:} Let $R=A[X,Y]$ be  a polynomial rings of indeterminates $X,Y;$ $\dim A=d>2.$  Then $R$ is  a finitely generated standard 2-graded algebra over $A$ with   $\deg X=(1,0),\deg Y=(0,1)$ and $$\dim D_{\frak m}(R)=\dim\biggl[\oplus_{n\ge 0}\dfrac{{\frak m}^n(XY)^n}{{\frak m}^{n+1}(XY)^n}\biggl]=\dim\biggl(\oplus_{n\ge 0}\dfrac{{\frak m}^n}{{\frak m}^{n+1}}\biggl)=\dim A.$$   
It can be verified that $X,Y$ is a weak-(FC)-sequence of $R$ with respect to ${\frak m}.$ Since $\dim D_{\frak m}(R/(X))= \dim (A/{\frak m})= 0$ and $d>2,$  $\dim D_{\frak m}(R/(X))<  \dim D_{\frak m}(R)-1.$ 
  
\vskip 0.2cm
From Theorem 3.3, in the case $s = 1$, we get the following result for a graded algebra $S = \bigoplus_{n \geq 0}S_n$.
\vskip 0.2cm
\noindent {\bf Corollary 3.6.} {\it Let $S = \bigoplus_{n \geq 0}S_n$ be a finitely generated standard graded algebra over $A$ such that $S_+ = \bigoplus_{n > 0}S_n$ is non-nilpotent and $J$ be an $\mathfrak{m}$-primary ideal of $A$.  Set $D_J(S) = \bigoplus_{n \geq 0}J^nS_n/J^{n + 1}S_n$ and $\dim D_\frak{m}(S) = \ell$. Suppose that $x_1, \ldots, x_q$ is a maximal weak-(FC)-sequence  in $S_1$ of $S$ with respect to $J$ satisfying the condition $\dim D_\frak{m}(S/(x_1, \ldots, x_q)S) = \ell - q$.  Then
\begin{enumerate}
\item[\rm (i)] $e(J, \ell - i - 1, i, S) \not= 0$ if and only if $i \leq q$ and $\dim A/((x_1, \ldots, x_i)S : S_1^\infty )_A = \ell - i$.

\item[\rm (ii)] If $e(J, \ell - i - 1, i, S) \not= 0$ then $ e(J, \ell - i - 1, i, S) = e_A(J, S_n/(x_1, \ldots, x_i)S_{n - 1})$ for all large $n$.
\end{enumerate} }
\vskip 0.2cm
 \noindent{\bfseries\itshape\textit Proof.} By Theorem 3.3(ii) we immediately get  (ii). Now let us to prove the part (i). The "if" part. Assume that $e(J, \ell - i - 1, i, S) \not= 0$. First, we  show that $i \leq q$. Assume the contrary that $i > q$. Since $x_1, \ldots, x_q$ is a weak-(FC)-sequence in $S_1$ of $S$ with respect to $J$, applying Proposition 3.2(i) by induction on $q$, 
$$0 \not= e(J, \ell - i - 1, i, S)  = e(J, \ell -i - 1, i - q, \bar{S}), $$
 where $\bar{S} = S/(x_1, \ldots, x_q)S$. Since  $e(J, \ell -i - 1, i - q, \bar{S}) \not= 0$ and $i - q > 0$, by Proposition 3.2(ii), there exists an element $x \in S_1$ such that $\bar{x}$ (the initial form of $x$ in $\bar{S}$) is a weak-(FC)-element of $\bar{S}$ with respect to $J$. By Proposition 3.2(i), 
$\dim D_\frak{m}(\bar{S}/x\bar{S}) = \ell - q - 1$.
 Hence $x_1, \ldots, x_q, x$ is a weak-(FC)-sequence in $S_1$ of $S$ with respect to $J$ satisfying the condition

\centerline {$\dim D_\frak{m}(S/(x_1, \ldots, x_q, x)S) = \ell - q - 1.$}
 \noindent We thus arrive at a contradiction. Hence $i \leq q$. Since $e(J, \ell - i - 1, i, S) \not= 0$, by Remark 3.4(iii), $ \dim A/((x_1, \ldots, x_i)S : S_1^\infty )_A = \ell - i.$
We turn to the proof of sufficiency. Suppose that $i \leq q$ and \\\centerline{$ \dim A/((x_1, \ldots, x_i)S : S_1^\infty )_A = \ell - i.$} Since $\dim D_\frak{m}(S/(x_1, \ldots, x_q)S) = \ell - q$, it follows, by Remark 3.4(i), that

\centerline{$\dim D_\frak{m}(S/(x_1, \ldots, x_i)S) = \ell - i$} 
\noindent for all $i \leq q$. 
 \noindent Since $x_1, \ldots, x_i$ is a weak-(FC)-sequence of $S$ with respect to $J$ satisfying the condition $$\dim D_\frak{m}(S/(x_1, \ldots, x_i)S) = \dim A/((x_1, \ldots, x_i)S : S_1^\infty )_A = \ell - i$$ by Theorem 3.3(i), $e(J, \ell - i - 1, i, S) \not= 0$. $\blacksquare$
\vskip 0.2cm
\noindent {\bf Example 3.7:} Let $R = A[X_1, X_2, \ldots, X_t]$ be the ring of polynomial in $t$
 indeterminates $X_1, X_2, \ldots, X_t$ with coefficients in $A$ ($\dim A = d > 0$). Then
$R =\bigoplus_{m \geq 0}R_m$ is a finitely generated standard graded algebra over $A$ (see Example 2.2). Let $J$ is an  $\mathfrak{m}$-primary ideal of $A$.  By Example 2.2,   $X_1, \ldots, X_t \in R_1$  is a weak-(FC)-sequence of $R$ with respect to $J$. Denote by $P(n, m)$ the polynomial of
$\ell_A(\frac{J^nR_m}{J^{n + 1}R_m})$.
We have
$$D_\mathfrak{m}(R) = \bigoplus_{T \geq 0} \frac{\mathfrak{m}^TR_T}{\mathfrak{m}^{T + 1}R_T} =  \frac{A[\mathfrak{m}X_1, \ldots, \mathfrak{m}X_t]}{\mathfrak{m}A[\mathfrak{m}X_1, \ldots, \mathfrak{m}X_t]} .$$
Since ht$\mathfrak{m} > 0$, $\dim D_\mathfrak{m}(R) =  \dim A + t - 1 = d + t - 1$. Hence  $\deg P(n, m) = d+ t - 2$.
It is clear that 
$R/(X_1, \ldots, X_i)R = A[X_{i + 1}, \ldots, X_t]$
for all $i \leq  t$. Hence $$\dim D_{\mathfrak{m}}(R/(X_1, \ldots, X_i)R) = \dim D_{\mathfrak{m}}(R) - i$$
 for all $i \leq t$.
 Let us calculate $e(J, k_0, k_1, R)$, with $k_0 + k_1 = d + t - 1$.
First, we consider  the case $k_1 \geq  t$. Since
$X_1, \ldots, X_{t - 1} $  is a weak-(FC)-sequence of $R$ with respect to $J$ and $\dim D_{\mathfrak{m}}(R/(X_1, \ldots, X_i)R) = \dim D_{\mathfrak{m}}(R) - i$ for all $i \leq t - 1$, by Remark 3.4(ii), 
$$e(J, k_0, k_1, R) = e(J, k_0, k_1 - (t - 1), R/(X_1, \ldots, X_{t - 1})R) = e(J, k_0, k_1 - t + 1, A[X_t]).$$
Denote by $Q(m, n)$ the polynomial of $\ell_A(\frac{J^nX_t^mA}{J^{n+1}X_t^mA})$. Since $X_t$ is regular element, $J^nX_t^mA \cong J^nA$. Thus, for all large $n, m$,
$$Q(n, m) = \ell_A(\frac{J^nX_t^mA}{J^{n+1}X_t^mA}) = \ell_A(\frac{J^nA}{J^{n+1}A}).$$
Hence $Q(n, m)$ is independent on $m$. Note that $e(J, k_0, k_1 - t + 1, A[X_t])$ is the coefficient of $\frac{1}{k_0!(k_1 - t + 1)!}n^{k_0}m^{k_1 - t + 1}$ in $Q(n, m)$. Since $k_1 - t + 1 > 0$, it follows that $$e(J, k_0, k_1, R) = e(J, k_0, k_1 - t + 1, A[X_t]) = 0.$$
In the case $k_1 <  t$, since $\dim D_{\mathfrak{m}}(R/(X_1, \ldots, X_{k_1})R) = \dim D_{\mathfrak{m}}(R) - k_1$, by Corollary 3.6(i), $e(J, k_0, k_1, R) \not= 0$ if and only if 
$$\dim A/((X_1, \ldots, X_{k_1})R: R_1^\infty )_A= d + t - 1 - k_1.$$
Since $X_1, \ldots, X_t$ are independent indeterminates, 
$$((X_1, \ldots, X_{k_1})R: R_1^\infty )_A \subset  ((X_1, \ldots, X_{k_1})R: ((X_{k_1 + 1}, \ldots, X_t)A)^\infty )_A = 0.$$
Hence $\dim A/((X_1, \ldots, X_{k_1})R: R_1^\infty )_A = \dim A = d$. Therefore,  
$e(J, k_0, k_1, R) \not= 0$ if and only if  $k_1 = t - 1$.  For $k_1 = t - 1$ (then $k_0 = d-1$), by Corollary 3.6(ii), we have
$$e(J, d-1, t - 1, R) = e_A(J, R_u/(X_1, \ldots, X_{t - 1})R_{u-1})$$
for all large $u$. Note that $R_u = (X_1, \ldots, X_t)^uA$. So $R_u/(X_1, \ldots, X_{t - 1})R_{u-1} = X_t^uA$. Thus $e(J, R_u/(X_1, \ldots, X_{t - 1})R_{u-1}) = e_A(J, X_t^uA)$.  Since $X_t^u$ is regular element in $A[X_t]$, $X_t^uA \cong A$.
Hence
$e(J, d-1, t - 1, R) = e_A(J, X_t^uA) = e_A(J, A).$
From the above facts we get 
\begin{equation*}
e(J, k_0, k_1, R) = \begin{cases} 
                            0  & \text{if }  k_1 \not= t-1\\
                            e_A(J, A) & \text{if }  k_1 = t-1
                            \end{cases}.
\end{equation*}
Therefore $$P(n, m) = \frac{e(J, A)}{(d-1)!(t - 1)!}n^{d-1}m^{t - 1} + \{\text{terms of lower degree}\}.$$

\vskip 0.2 cm
\noindent{\bf Remark 3.8.} Example 3.5 and Example 3.7 showed that for any weak-(FC)-sequence $x_1, \ldots, x_t$ of $S$  with respect  to $J,$
one can get 
$$\dim D_{\frak m}(S/(x_1, \ldots, x_{t})) < \dim D_{\frak m}(S)-t,$$ and 
$\dim A/((x_1, \ldots, x_{t})S: S_1^\infty )_A = \dim [A/(0:S_1^\infty )_A] -t$ although $$\dim A/((x_1, \ldots, x_{i})S: S_1^\infty )_A \not= \dim [A/(0:S_1^\infty )_A] -i$$ for some $i< t.$ That is a difference of weak-(FC)-sequences in  graded algebras and weak-(FC)-sequences of ideals in local rings.      

\vskip 0.2cm
\noindent {\bf 4. Applications}

As an application of Theorem 3.3, this section devoted to the discussion of  mixed multiplicities of arbitrary ideals in local rings.

Throughout this section, let $(A, \frak{m})$ denote a Noetherian local ring with maximal ideal 
$\mathfrak{m}$, infinite residue $k = A/\mathfrak{m},$
  and  an ideal  $\frak m$-primary  $J$, and $I_1,\ldots, I_s$   ideals of $A$ such that $I = I_1\cdots I_s$ is  non-nilpotent.  Set $S = A[I_1t_1,\ldots,I_st_s]$. Then   
$$D_J(S) = \bigoplus_{n \geq  0} \frac{(JI)^n}{J(JI)^n} \;\; \text {and} $$ 
$$ \ell_A(\frac{J^{n_0}S_{(n_1, \ldots, n_s)}}{J^{n_0 + 1}S_{(n_1, \ldots, n_s)}}) = \ell_A(\frac{J^{n_0}I_1^{n_1}\cdots I_s^{n_s}}{J^{n_0 + 1}I_1^{n_1}\cdots I_s^{n_s}})$$
is a polynomial of total degree $\dim D_J(S) -1$  for all large $n_0, n_1, \ldots, n_s.$ 
 By  Proposition 3.1 in [Vi],  the  degree  of  this polynomial  is $\dim A/0:I^\infty -1$. Hence $\dim D_J(S) = \dim A/0:I^\infty.$
Set $\dim A/0:I^\infty = \ell.$ 
In this case, $e(J, k_0, k_1, \ldots, k_s, S)$ for $ k_0 + k_1 + \cdots + k_s = \ell - 1$ is called the  mixed multiplicity of ideals $(J,I_1,\ldots, I_s)$ of type $(k_0, k_1, \ldots, k_s)$ and one put  $$e(J, k_0, k_1, \ldots, k_s, S)= e(J^{[k_0+1]}, I_1^{[k_1]},\ldots, I_s^{[k_s]}, A)$$ (see [Ve2] or [HHRT]). By using the concept of (FC)-sequences of ideals, one 
 transmuted  mixed multiplicities of a set of arbitrary ideals into
 Hilbert-Samuel multiplicities [Vi].    
\vskip 0.2cm
\noindent {\bf Definition 4.1} [see Definition, Vi]{\bf.}  {\it Let  $I_1,\ldots, I_s$  be  ideals such that $I = I_1\cdots I_s$ is a non nilpotent ideal.
A element $x \in A$ is called an (FC)-element of $A$ with respect to $(I_1,\ldots, I_s)$ 
if there exists $i \in \{1, 2, \ldots, s\}$ such that $x \in I_i$ and
\begin{list}{}{ \setlength{\leftmargin}{1.8cm}\setlength{\labelwidth}{1.3cm} }
\item [\rm (FC$_1$): ] $(x)\cap I_1^{n_1}\cdots I_i^{n_i}\cdots I_s^{n_s} 
= xI_1^{n_1}\cdots I_i^{n_i-1}\cdots I_s^{n_s}$
for all large $n_1, \ldots, n_s$.
\item [\rm (FC$_2$): ] $x$ is a filter-regular element with respect to $I,$ i.e., $0:x\subseteq 0:I^\infty.$

 \item [\rm (FC$_3$): ] $\dim A/[(x):I^\infty]=\dim A/0:I^\infty-1.$
\end{list} 

We call $x$ a {\it weak-(FC)-element} with respect to $(I_1,\ldots, I_s)$ if $x$ satisfies conditions $(FC_1)$ and $(FC_2).$

Let $x_1, \ldots, x_t$ be a sequence in $A$. For each $i = 0, 1, \ldots, t - 1 $, set $A_i = A/(x_1, \ldots, x_{i})S$, $\bar{I}_j = I_j[A/(x_1, \ldots, x_{i})]$, $\bar{x}_{i + 1}$ the image of $x_{i + 1}$ in $A_i$. Then 

$x_1, \ldots, x_t$ is called a weak-(FC)-sequence  of $A$ with respect to $(I_1,\ldots, I_s)$ if $\bar{x}_{i + 1}$ is a weak-(FC)-element  of $A_i$ with respect to $(\bar{I}_1,\ldots, \bar{I}_s)$ for all $i = 0, 1, \ldots, t - 1$.

$x_1, \ldots, x_t$ is called an (FC)-sequence  of $A$ with respect to $(I_1,\ldots, I_s)$ if $\bar{x}_{i + 1}$ is an (FC)-element  of $A_i$ with respect to $(\bar{I}_1,\ldots, \bar{I}_s)$ for all $i = 0, 1, \ldots, t - 1$.

A weak-(FC)-sequence $x_1,\ldots,x_t$  is called a {\it maximal} weak-(FC)-sequence if $IA_{t-1}$ is a non-nilpotent ideal of $A_{t-1}$ and $IA_t$ is a nilpotent ideal of $A_t.$}
\vskip 0.2 cm
\noindent {\bf Remark 4.2.} 
\begin{itemize}
\item[(i)] The  condition (FC$_1$) in Definition 4.1 is a weaker condition than the condition (FC$_1$) of definition of (FC)-element in [Vi]. 

\item[(ii)]  If $x\in I_i$ is a weak-(FC)-element with respect to $(J, I_1,\ldots, I_s),$ then it can be verified that $x$ is also a weak-(FC)-element of $S$ with respect to $J$ as in  Definition 2.1. 

\item[(iii)] If $x_1, \ldots, x_t$ is an (FC)-sequence with respect to $(J, I_1,\ldots, I_s),$ then from the condition (FC$_3$) we follow that $\dim A/((x_1, \ldots, x_t)S: {S}_{(1, \ldots, 1)}^\infty )_A= \ell - t.$ Hence  $$\dim D_J(S/(x_1, \ldots, x_t)S) = \dim A/((x_1, \ldots, x_t)S: {S}_{(1, \ldots, 1)}^\infty )_A= \ell - t$$ that as in the state of Theorem 3.3(i).
\item[(iv)] By Lemma 3.1, $e(J, k_0, 0, \ldots, 0, S) \not= 0 $ if and only if
$\dim A/(0 : S_{(1, \ldots, 1)}^\infty)_A = \ell.$
 In this case, $e(J, k_0, 0, \ldots, 0, S) = e_A(J, S_{(n, \ldots, n)})$ for all large $n$. But since $\dim A/(0 : S_{(1, \ldots, 1)}^\infty)_A = \dim A/0:I^\infty,$ $\dim A/(0 : S_{(1, \ldots, 1)}^\infty)_A = \ell.$ Hence 
$e(J, k_0, 0, \ldots, 0, S) = e_A(J, S_{(n, \ldots, n)})$  for all large $n.$ It is a plain fact that $e_A(J, S_{(n, \ldots, n)}) = e_A(J, I^n).$  
On the other hand by the proof of Lemma 3.2 [Vi], $e_A(J, I^n)= e_A(J,A/0:I^\infty)$ for all large $n.$ Hence $e(J, k_0, 0, \ldots, 0, S) = e_A(J,A/0:I^\infty).$   

\end{itemize}
 
 Then as an immediate consequence of  Theorem 3.3, we obtained a more favorite  result than [Theorem 3.4, Vi](see Remark 4.2 (i)) as follows.
 \vskip 0.2cm
\noindent {\bf Theorem 4.3} [see Theorem 3.4, Vi]{\bf.} {\it Let $(A, \frak{m})$ denote a Noetherian local ring with maximal ideal 
$\mathfrak{m}$, infinite residue $k = A/\mathfrak{m},$
  and  an ideal  $\frak m$-primary  $J$, and $I_1,\ldots, I_s$   ideals of $A$ such that $I = I_1\cdots I_s$ is non nilpotent. Then the following statements hold.
\begin{enumerate}
\item[\rm (i)] $e(J^{[k_0+1]}, I_1^{[k_1]},\ldots, I_s^{[k_s]}, A) \not= 0$ if and only if there exists a weak-(FC)-sequence 
$x_1,  \ldots, x_t$  with respect to $(J, I_1,\ldots, I_s)$ 
consisting of $k_1$ elements of $I_1$, ..., $k_s$ elements of $I_s$ and 
$\dim A/(x_1,  \ldots, x_t) : I^\infty = \dim A/0: I^\infty -t.$
\item[\rm (ii)] Suppose that $e(J^{[k_0+1]}, I_1^{[k_1]},\ldots, I_s^{[k_s]}, A) \not= 0$ and $x_1,  \ldots, x_t$ 
 is a weak-(FC)-sequence  with respect to $(J, I_1,\ldots, I_s)$ consisting of $k_1$ elements of $I_1$, ..., $k_s$ elements of $I_s$. Set $\bar{A} = A/(x_1, \ldots, x_t):I^\infty$. Then
$$e(J^{[k_0+1]}, I_1^{[k_1]},\ldots, I_s^{[k_s]}, A) = e_A(J,\bar{A}).$$
\end{enumerate}}

Note that one can  get this result by a minor improvement in the proof of [Proposition 3.3, Vi](see [DV1]). Moreover,  the filtration version of Theorem 4.3 is proved in [DV2].
   
Recently, [DV1] and [DV2] showed that from  Theorem 4.3 one rediscover the earlier  result of Trung and Verma [TV] on mixed multiplicities of ideals. This fact proved  that Theorem 3.3 covers the main results in [Vi] and [TV].

\end{document}